\documentclass{amsart}

\usepackage{amssymb}

\newtheorem{theorem}{Theorem}[section]
\newtheorem{lemma}[theorem]{Lemma}
\newtheorem{proposition}[theorem]{Proposition}
\newtheorem{corollary}[theorem]{Corollary}

\theoremstyle{definition}
\newtheorem{definition}[theorem]{Definition}

\theoremstyle{remark}
\newtheorem{remark}[theorem]{Remark}

\numberwithin{equation}{section}

\newcommand{\h}{H(X)}
\newcommand{\f}{\varphi}
\newcommand{\M}{\mathcal M(\varphi)}
\newcommand{\e}{\varepsilon}

\begin{document}

\title[Full groups of Cantor minimal systems]{Full groups, flip conjugacy,
and orbit equivalence of  Cantor minimal systems}


\author{S. Bezuglyi}
\address{Insitute for Low Temperature Physics, National Academy of Sciences
of Ukraine, Kharkov, Ukraine}
\email{bezuglyi@ilt.kharkov.ua}

\author{K. Medynets}
\address{Insitute for Low Temperature Physics, National Academy of Sciences
of Ukraine, Kharkov, Ukraine}
 \email{medynets@ilt.kharkov.ua}
\thanks{The second named author is supported by INTAS YSF-05-109-5315}

\subjclass[2000]{Primary 37B05, Secondary 20B99}

\date{}
\keywords{Cantor set, minimal homeomorphism, full groups, commutator,
involution}



\begin{abstract} In the paper, we consider the full group $[\varphi]$
and topological full group $[[\varphi]]$ of a Cantor minimal system
$(X,\f)$. We prove that the commutator subgroups $D([\f])$ and $D([[\f]])$
are  simple and show that the groups $D([\f])$ and $D([[\f]])$ completely
determine the class of orbit equivalence and flip conjugacy of $\f$,
respectively. These results improve the classification found in
\cite{gps:1999}.
As a corollary of the technique used, we establish the fact
that $\f$ can be written as a product of three involutions from $[\f]$.
\end{abstract}

\maketitle


\section{Introduction}

One of the most remarkable results of ergodic theory is Dye's theorem
which states that any two ergodic finite measure-preserving automorphisms
of a Lebesgue space are orbit equivalent  and, as a corollary,
their full groups are isomorphic \cite{Dye1}. Dye also proved that the
full group is a
complete invariant of orbit equivalence for ergodic finite
measure-preserving actions of countable groups  \cite{dye}. Later, full
groups have been studied in numerous papers from different points of view.
In particular, the analogues of the theorems of Dye were established for
infinite measure-preserving and non-singular automorphisms of a standard
measure space.

The ideas developed in ergodic theory have been successfully applied to the
study of orbit equivalence in Cantor and Borel dynamics. Giordano, Putnam,
and Skau considered the notions of the full group $[\f]$ and topological
full group $[[\f]]$ of a Cantor minimal system $(X,\f)$ and showed that
these groups completely determine the classes of orbit equivalence and flip
conjugacy of $\f$ \cite{gps:1999}. In other words, they proved that, for
minimal homeomorphisms $\f_1$ and $\f_2$, any algebraic isomorphism of full
groups $[\f_1]$ and $[\f_2]$ is spatially generated. Recently, Miller and
Rosendal have shown that Dye's theorem holds in the context of Borel
dynamics: two Borel aperiodic actions of countable groups are orbit
equivalent if and only if their full groups are isomorphic
\cite{MillerRosendal}. We should also mention that there are several papers
where the algebraic and topological structures of full groups of ergodic
automorphisms have been studied. In particular, Eigen proved that the full
group of an ergodic finite measure-preserving automorphism is simple, i.e.
it has no proper normal subgroups \cite{eigen:1981}. Notice that it is
still an open problem to show the simplicity of the full group $[\f]$
generated by a minimal homeomorphism of a Cantor set. On the other hand, it
is known that the topological full group $[[\f]]$ has a proper normal
subgroup and so is not simple.

The goal of this paper is to prove that the class of orbit equivalence and
flip conjugacy of a Cantor minimal system can be determined by proper
simple subgroups of the full group and topological full group.
To do this,  we focus our study on the commutator subgroups $D([\f])$
and $D([[\f]])$ of $[\f]$ and $[[\f]]$, respectively. First of all, we prove
in Section \ref{Section_Commutator_Subgroups} that the groups $D([\f])$
and $D([[\f]])$ are simple (Theorem \ref{simplicity}). Then we show that
the commutator subgroup $D([[\f]])$ is a complete invariant of flip
conjugacy (Theorem \ref{Theorem_FlipConjugacy}) and the group $D([\f])$ is a
complete invariant of orbit equivalence (Theorem
\ref{Theorem_OrbitEquivalence}). In such a way, these results make more
precise the characterization of orbit equivalence and flip conjugacy found
in \cite{gps:1999}). In particular, Theorem \ref{Theorem_FlipConjugacy}
contains a new proof of the fact that $[[\f]]$ is a complete invariant of
flip conjugacy.  To show that the group $D([\f])$ determines the class of
orbit equivalence, we follow the idea of the proof of Corollary 4.6 from
\cite{gps:1999}. The key point of our approach is the fact that every
involution with clopen support belongs to $D([\f])$ (Corollary
\ref{ClopenInvol}). In its turn, this result is based on the fact that
every homeomorphism from $[\f]$, which is minimal on its support, can be
written as a product of five commutators from $[\f]$ (Theorem
\ref{ProductCommutators}).

In  Section \ref{Section_ProductInvolutions}, we also answer the question
of representation of a minimal homeomorphism $\f$ as a product of
involutions from $[\f]$. We note that this problem would be trivial if we
knew that $[\f]$ is a simple group. The problem of writing every element
of a transformation group as a product of involutions has been considered in
measurable, Borel, Boolean dynamics.
Apparently, the first results appeared in the paper of Anderson
\cite{anderson:1958}, where it was shown that every homeomorphism of a
Cantor set is a product of six involutions from $\h$, the group of all
homeomorphisms of a Cantor set $X$. The technique used by Anderson also
works for the group of non-singular transformations of a Lebesgue space
(this fact was mentioned in \cite{eigen:1981}). In \cite{fathi}, Fathi
suggested a new approach to this problem which allowed him to show that
the group $Aut(Y,\mu)$ of finite measure-preserving automorphisms of a
Lebesgue space $(Y,\mu)$ is  simple, and that every $f\in Aut(Y,\mu)$ is a
product of 10 involutions. His ideas were then used in \cite{eigen:1981}
to obtain a similar result for the full group of an ergodic automorphism
of a Lebesgue space. Later, Ryzhikov improved Fathi's theorem and showed
that every ergodic automorphism $f\in Aut(X,\mu)$ is, in fact, a product
of three involutions from its full group \cite{ryzjikov:1985}. We notice that
the number `three' is the least possible for such a representation by
involutions taken from the full group. Miller showed that an automorphism
$f$ is the product of two involutions from its full group if and only if
$f$ is dissipative \cite{miller:thesis}. We mention also the works
\cite{fremlin}, \cite{ryzhikov:1993} and \cite{choksi_prasad} where these
questions were studied for automorphisms of complete Boolean algebras and
homogeneous measure algebras. Developing ideas of Fathi and Ryzhikov, we
prove that every minimal homeomorphism $\f$ of a Cantor set is the product
of three involutions from $[\f]$ (Theorem
\ref{Theorem_Product3involutions}).

\section{Preliminaries}\label{Section-Preliminaries}
Let $X$ denote a Cantor set, i.e. a 0-dimensional compact metric space
without isolated points. Denote by $H(X)$ the group of all homeomorphisms
of $X$. For a homeomorphism $\f\in H(X)$ and a point $x\in X$, let
$Orb_\f(x)=\{\f^n(x) : n\in\mathbb Z\}$ denote the $\f$-orbit of $x$. The
open set $supp(\f)=\{x\in X : \f(x)\neq x\}$ is called the {\it support of}
$\f$.

A homeomorphism $\varphi\in H(X)$ is called {\it periodic} if for each
$x\in X$ there exists $n>0$ such that $\f^n(x)=x$, and  if  the
$\f$-orbit of $x$  is infinite for all $x$, then $\f$ is called {\it aperiodic}.
A homeomorphism $\f$ is called {\it  minimal} if the $\varphi$-orbit of
every point is dense in $X$.  Let
$[\varphi]$ denote the set of all homeomorphisms $f\in H(X)$ such
that $f\in Orb_\varphi(x)$ for all $x\in X$. The set $[\f]$ is
called the {\it full group of $\varphi$}. If $f\in [\f]$, then there is a
function $n_f:X\rightarrow \mathbb Z$ such that $f(x)=\f^{n_f(x)}(x)$ for
all $x\in X$, which is called the {\it cocycle associated to} $f$. The
{\it topological full group $[[\f]]$ of $\f$} is defined as the set of all
$f\in [\f]$ such that $n_f$ is continuous.  We refer the reader to
\cite{gps:1999} and \cite{glasner_weiss:1995} for an in-depth study of
various properties of full groups associated to a minimal homeomorphism of
a Cantor set.

Given a group $G$, denote by $D(G)$ the subgroup generated by all elements
of the form $[f,g]:=fgf^{-1}g^{-1}$, $f,g\in G$. The group $D(G)$ is
called the {\it commutator subgroup of}  $G$.

Let $\f\in\h$ and let $A$ be a clopen set. A point $x\in A$ is called {\it
recurrent} for $\f$ if  there is $n>0$ such that $\f^n(x)\in A$. Observe that if a
clopen set $A$ consists of recurrent points, then  the function $n_A:
A\rightarrow \mathbb N$ given by $n_A(x)=\min\{n>0 : \f^n(x)\in A\}$ is
well-defined and continuous. Then $n_A$ is called the {\it the  first
return function}. Suppose that a clopen set $A$ meets every $\f$-orbit and
consists of recurrent points. Set $A_k=\{x\in A: n_A(x)=k\}$. Hence, $X$
can be decomposed into clopen sets
$$X=\bigcup_{k\geq 1}\bigcup_{i=0}^{k-1}\f^{i}(A_k).$$
This decomposition is called a {\it Kakutani-Rokhlin (K-R) partition}
built over the set $A$. By a {\it $\f$-tower}, we mean a non-empty
disjoint family $\xi=(A_k,\f A_k,\ldots,\f^{k-1}A_k)$. The set $A_k$ is
called the {\it base} of $\xi$ and $k$ is called the {\it height} of
$\xi$.

The following simple proposition explains the construction of the induced
map.
\begin{proposition}\label{InducedMap} Let $\f\in\h$ and let a clopen set $A$
meet every $\f$-orbit. Then:

 (1) $A$ consists of recurrent points;

(2) the homeomorphism $\f_A\in \h$ given by $\f_A(x)=\f^{n_A(x)}(x)$ for
$x\in A$ and $\f_A(x)=x$ for $x\in X\setminus A$ belongs to $[[\f]]$;

(3) the homeomorphism $g=\f_A^{-1}\f$ is periodic.
\end{proposition}
{\it Proof.} (1) Since $A$ meets every $\f$-orbit, we have that
$X=\bigcup_{i\in\mathbb Z}\f^i(A)$. By compactness of $X$, we get that
$X=\bigcup_{i=0}^m\f^i(A)$ for some $m$. This implies that $A$ consists of
recurrent points. The statement (2) is trivial and the proof of (3) is
straightforward. \hfill$\square$

\smallskip
For a homeomorphism $\f$, let us denote by $\mathcal M(\varphi)$ the set
of all Borel probability $\varphi$-invariant measures on $X$. Clearly, if
$\gamma\in [\varphi]$, then $\gamma\circ\mu=\mu$ for every $\mu\in
\mathcal M(\varphi)$. The following theorem answers the question when two
clopen sets can be mapped onto each other by an element of the full group.
This result will be one of our main tools in the study of full groups.

\begin{theorem}[\cite{glasner_weiss:1995}]\label{equivalent-clopen-sets}
Let $(X,\varphi)$ be a Cantor minimal system.

(1) If $A,B$ are clopen subsets of $X$ such that $\mu(B)<\mu(A)$ for every
$\mu\in \M$, then there exists $\alpha\in [[\f]]$ with $\alpha(B)\subset
A$. Moreover, $\alpha$ can be chosen such that $\alpha^2=id$ and $\alpha
|(B\cup \alpha(B))^c=id$ where $E^c := X\setminus E$.

(2) If $A,B$ are clopen sets with $\mu(A)=\mu(B)$ for all $\mu\in \M$,
then there exits $\alpha\in [\f]$ such that $\alpha(B)=A$. Moreover,
$\alpha$ can be chosen such that $\alpha^2=id$, $\alpha$ has clopen support,
and $\alpha |(B\cup \alpha(B))^c=id$.
\end{theorem}

We will need a generalization of a theorem from \cite{akin:1999}.

\begin{proposition}\label{propertiesOfInvariantMeasures} Let $(X,\f)$ be a
Cantor minimal system and let $d$ be a metric on $X$ compatible with the topology.

(1) For any $\e> 0$ there exists $\delta>0$ such that  if the $d$-diameter
of a clopen set $A$ is less than $\delta$, then $\mu(A)<\e$ for every
$\mu\in\M$.

(2) If $A$ is a clopen set, then $inf\{\mu(A) : \mu\in \M\}>0$.
\end{proposition}

{\it Proof}. (1) Assume the converse: there exists $\e_0> 0$ such that for
every $n$ there exist a clopen set $A_n$ with $diam(A_n)<1/n$ and a
measure $\mu_n\in M(\f)$ such that $\mu_n(A_n)\geq \e_0$. By compactness
of $M(\f)$ and $X$, we may assume, without loss of generality, that
$\mu_n\to\mu\in M(\f)$ and there exists $x_0\in X$ such that every
neighborhood of $x_0$ contains all but a finite number of sets $A_n$.

Consider any clopen neighborhood $O$ of $x_0$. Then $A_n\subseteq O$ and
$\mu_n(O)\geq\mu_n(A_n)\geq\e_0$ for sufficiently large  $n$.  As
$\mu_n\to\mu$, we get that $\mu(O)\geq \e_0$. Hence $\mu(\{x_0\})>0$, a
contradiction.

(2) It easily follows from the minimality of $\f$ (see also
\cite{glasner_weiss:1995}). \hfill$\square$
\smallskip

From now on,  the pair $(X,\f)$  will always denote a Cantor minimal
system.
%
%

\section{Commutator subgroups}\label{Section_Commutator_Subgroups}
In this section, we show that, for a minimal homeomorphism $\f$,
 the commutator subgroups $D([\f])$ and
$D([[\f]])$ are  simple. In our proofs, we follow the arguments
used in \cite{fathi} and \cite{eigen:1981}.

The following simple statement describes the properties of periodic
homeomorphisms. More detailed descriptions of  periodic homeomorphisms
from $[[\f]]$ and $[\f]$ can be found in \cite{bezuglyi_dooley_kwiatkowski}.

\begin{lemma}\label{PeriodicProperties} Let
$f\in \h$ and $X_n:= \{x\in X : |Orb_f(x)|=n\}$ for $1\leq n<\infty$.

(1) If $f\in [[\f]]$ with $\f$ a minimal homeomorphism, then $X_n$ is clopen.

(2) If $X_n$ is clopen, then there exists a clopen set $X_n^0$ such that
$X_n=\bigcup_{i=0}^{n-1}f^i(X_n^0)$ a disjoint union.
\end{lemma}

{\it Proof.} (1) Notice that $X_{\leq n}: =\{x\in X : |Orb_f(x)|\leq n\}$
is closed. Since $f\in [[\f]]$, the associated cocycle $n_f$ is
continuous. This implies that the set $X_n$ is open for each $1\leq  n<\infty$.
Hence, $X_n$ is clopen.

(2) See Lemma 3.2 in \cite{bezuglyi_dooley_kwiatkowski}. \hfill$\square$
\medskip

The following lemma states that every homeomorphism from the full group
can be written as a product of homeomorphisms which have `small' supports.

\begin{lemma}\label{ProductSmallGenerators} Let $\Gamma$ denote either $[\f]$ or
$[[\f]]$. Given $\delta>0$ and $g\in \Gamma$, there exist
$g_1,\ldots,g_m\in \Gamma$ and clopen sets $E_1,\ldots,E_m$ such that

(1) $g=g_1,\ldots,g_m$;

(2) $supp(g_i)\subseteq E_i$ and $\mu(E_i)<\delta$ for all $\mu\in \M$ and
$i=1,\ldots,m$.
\end{lemma}

{\it Proof.} Suppose first that $\Gamma=[\f]$. Using Proposition
\ref{propertiesOfInvariantMeasures}, take any decomposition of $X$ into
disjoint clopen sets $X=B_1\sqcup\ldots\sqcup B_m$ such that
$\mu(B_i)<\delta/2$ for any $\mu\in \M$ and $i=1,2,\ldots,m$.

Given $f\in [\f]$, we have that  $\mu(B\setminus f(B))=\mu (f(B)\setminus
B)$ for every $\mu\in \M$ and every clopen set $B$. Therefore, by
Theorem \ref{equivalent-clopen-sets}, there exists $g_1\in [\f]$ such that
$g_1|B_1= g|B_1$, $g_1(g(B_1)\setminus B_1)=B_1\setminus g(B_1)$, and
$supp(g_1)\subseteq B_1\cup g(B_1)$. Setting $E_1=B_1\cup g(B_1)$, we obtain
that $\mu(supp(g_1))\leq\mu(E_1)<\delta$ for all $\mu\in \M$.

Let $g_1'=g_1^{-1}g$. Clearly, $supp(g'_1)\subseteq B_2\sqcup\ldots \sqcup
B_m$. Find $g_2\in [\f]$ such that $g_2|B_2=g_1'|B_2$ and
$supp(g_2)\subseteq B_2\cup g_1'(B_2)$. Let $E_2=B_2\cup g_1'(B_2)$. Then $\mu(E_2)<\delta$ for all $\mu\in \M$.
Next, consider $g_ 2':=g_2^{-1}g_1'$. Clearly, $supp(g_2')\subseteq
B_3\sqcup\ldots\sqcup B_m$ and $g=g_1g_2g_2'$. Repeating the above
argument for each set $B_i$, we construct a family of homeomorphisms
$g_i\in [\f]$,\ $i=1,\ldots,m$, such that $g=g_1\ldots g_m$ and
$supp(g_i)\subseteq B_i\cup g_{i-1}'(B_i)$. Here $g_m=g_{m-1}'=g_{m-1}^{-1}g_{m-2}'$. Setting $E_i=B_i\cup
g_{i-1}'(B_i)$, we complete the proof for the case $\Gamma = [\f]$.

Suppose now that $\Gamma=[[\f]]$. If $g$ is periodic, then by Lemma
\ref{PeriodicProperties} and compactness of $X$, we can  decompose $X$
into a finite number of clopen sets
$$
X=\bigcup_{i\in I}\bigcup_{j=0}^{n_i-1}g^j(X_{i}^0),
$$
where $g^{n_i}(x)=x$ for each $x\in X_i^0$. By Proposition
\ref{propertiesOfInvariantMeasures}, we can divide each set $X_{i}^0$ into
clopen sets $\{A_{i,1}^0\ldots,A_{i,m_i}^0\}$ such that $\mu(\mathcal
A_{i,j})<\delta$ for all $\mu\in  \M$, where $\mathcal
A_{i,j}=A_{i,j}^0\cup\ldots\cup g^{n_i-1}A_{i,j}^0$. Set $g_{i,j}x=gx$
whenever $x\in \mathcal A_{i,j}$ and $g_{i,j}x=x$ elsewhere. Clearly,
every $g_{i,j}\in [[\f]]$ and $g$ is the product of the commuting elements
$g_{i,j}$.

Consider the case when $g$ is not periodic. Choose an integer $k>0$ such
that $1/k<\delta$. As $g\in [[\f]]$, the set $X_{\geq k}:= \{x\in X :
|Orb_g(x)|\geq k\}$ is clopen. Using the arguments of \cite[Proposition
3]{bezuglyi_dooley_medynets} (see also \cite[Lemma 2.2]{medynets:2006}),
we can show that there exists a clopen set $B\subset X_{\geq k}$ such that
$g^i(B)\cap B=\emptyset$ for $i=0,\ldots,k-1$ and $B$ meets every
$g|X_{\geq k}$-orbit. It follows that $\mu(B)\leq 1/k<\delta$ for all
$\mu\in \M$. By our choice of the set $B$, the induced
homeomorphism $g_B$, defined as in Proposition \ref{InducedMap}, belongs to
$[[\f]]$. Moreover, $\mu(supp(g_B))<\delta$ for all $\mu\in \M$. Observe
that $g_1=g_B^{-1}g$ is periodic.  We can use the decomposition obtained
for periodic homeomorphism $g_1$ to complete the proof of the lemma.
\hfill$\square$

\begin{lemma}\label{generators} Suppose that  $G$ is a group and $H$ is a normal
subgroup of $G$. If $g_1,\ldots,g_n,h_1,\ldots,h_m$ are such elements from $G$
that $[g_i,h_j], [g_i,g_j]$, and $[h_i,h_j]$ belong to $H$ for any $i,j$, then
the element $[g_1\cdots g_n,h_1\cdots h_m]$ also belongs to $H$.
\end{lemma}
{\it Proof.} Note that
$$[g_1g_2,h_i]=g_1[g_2,h_i]g_1^{-1}[g_1,h_i],
$$
$$
[g_j,h_1h_2]=[g_j,h_1]h_1[g_j,h_2]h_1^{-1}.
$$

As $H$ is a normal subgroup, then $[g_1g_2,h_i]$ and $[g_i,h_1h_2]$ belong
to $H$. Hence $[g_1g_2,h_1h_2]\in H$. The proof can be completed by
induction. We leave the details to the reader. \hfill$\square$
\medskip

Now we are ready to show that the commutator subgroups $D([\f])$
and $D([[\f]])$ are simple. We should also mention that the
simplicity of $D([[\f]])$  was first established by Matui
\cite[Theorem 4.9.]{matui}, but with a completely different
technique.

\begin{theorem}\label{simplicity} Let $(X,\f)$ be a Cantor minimal
system.

(1) If $H$ is a normal subgroup of $[\f]$ (of $D([\f])$), then
$D([\f])\subseteq H$.

(2) If $H$ is normal subgroup of $[[\f]]$ (of $D([[\f]])$), then
$D([[\f]])\subseteq H$.

In particluar, the groups $D([\f])$ and $D([[\f]])$ are simple.
\end{theorem}

{\it Proof.} First of all notice that if $2\mu(B)<\mu(A)$ for any $\mu\in M(\f)$ with $A$ and $B$ clopen sets, then there exists $\alpha \in D([[\f]])$ such that $\alpha(B)\subset A$. Indeed, by setting $\alpha=id$ on $A\cap B$, we may assume that $A\cap B=\emptyset$. Applying  Theorem \ref{equivalent-clopen-sets} twice, we find two involutions $\alpha_1,\alpha_2\in [[\f]]$ such that $\alpha_1(B)\subset A$, $\alpha_2(\alpha_1(B))\subset A\setminus \alpha_1(B)$, and $supp(\alpha_1)=B\cup \alpha_1(B)$,
 $supp(\alpha_2)=\alpha_1(B)\cup \alpha_2(\alpha_1(B))$.
  Set $\alpha=\alpha_1\alpha_2$.
  Therefore, $\alpha(B)=\alpha_1(B)\subset A$. Since $\alpha_2=\alpha \alpha_1^{-1} \alpha^{-1}$, we get that $\alpha=\alpha_1\alpha_2=[\alpha_1,\alpha]$.

We must show that $[g,h]\in H$ for any $g,h\in \Gamma$, where $\Gamma$
denotes one of the groups $D([[\f]])$, $[[\f]]$, $D([\f])$, or $[\f]$. Take
any non-trivial element $f\in H$. Then there exists a clopen set
$E\subseteq X$ such that $f(E)\cap E=\emptyset$. Proposition
\ref{propertiesOfInvariantMeasures} implies that
$\delta=\frac{1}{2}\inf\{\mu(E) : \mu\in \M \}>0$.

By Lemma \ref{ProductSmallGenerators}, write $g$ and $h$ as
$g=g_1,\ldots,g_n$ and $h=h_1,\ldots,h_m$, such that (1)
$g_i,h_j\in\Gamma$ and (2) there exist clopen sets $E_i(g)$ and $E_j(h)$
with $supp(g_i)\subseteq E_i(g)$, $supp(h_j)\subseteq E_j(h)$ and
$\mu(E_i(g)\cup E_j(h))< \delta$, for every $\mu\in \M$ and every $i,j$.
Due to Lemma \ref{generators}, it is sufficient to prove that
$[g_i,h_j]\in H$.

For convenience we omit subindexes $i$ and $j$. Consider any homeomorphisms
$g,h\in \Gamma$ such that $supp(g)\cup supp(h)\subseteq F$, where $F$ is a
clopen set with $\mu(F)<\delta$ for all $\mu\in \M$. As above, find a
homeomorphism $\alpha\in D([[\f]])\subset\Gamma$ with $\alpha(F)\subseteq
E$. Since $H$ is a normal subgroup of $\Gamma$, the element $q=\alpha^{-1}
f\alpha\in H$.

Since $H$ is a normal subgroup, we have that $\widehat{h}=[h,q]=hqh^{-1}q^{-1}\in H$.
Analogously,  $[g,\widehat h]\in H$. Since $q(F)\cap F=\emptyset$, we see
that $g^{-1}$ and $ qh^{-1}q^{-1}$ commute. Then
$$[g,\widehat h]= g (hqh^{-1}q^{-1})
g^{-1}(qhq^{-1}h^{-1})= ghg^{-1} qh^{-1}q^{-1}qhq^{-1}h^{-1}=[g,h]\in H.
$$
This completes the proof. \hfill$\square$

\begin{remark}\label{Remark_NormalSubgroups} For a measure $\mu\in \M$,
set $[[\f]]_0=\{\gamma\in [[\f]] : \int_Xn_\gamma d\mu=0\}$. The
definition of $[[\f]]_0$ does not depend on the choice of $\mu$
\cite[Section 5 ]{gps:1999}. As proved in \cite{gps:1999}, the
group $[[\f]]_0$ completely determines the class of flip conjugacy
of $\f$. Clearly, $[[\f]]_0$ is a proper normal subgroup of
$[[\f]]$. Therefore, by Theorem \ref{simplicity},
$D([[\f]])\subseteq [[\f]]_0$ and $D([[\f]])=D([[\f]]_0)$. In
\cite{gps:1999} the authors asked if $[[\f]]_0$ is a simple group.
This would show that the class of flip conjugacy is determined by
a countable simple group. However, $[[\f]]_0$ is not simple, in
general. Matui proved that the simplicity of $[[\f]]_0$ is
equivalent to the 2-divisibility of the dimension group
$K^0(X,\f)$ \cite{matui}. Nevertheless, in Section
\ref{SectionFlipConjugacy} we will show that $D([[\f]])$ is a
complete invariant for flip conjugacy.
\end{remark}

%
%

\section{Product of involutions}\label{Section_ProductInvolutions}
In the section, we show that a minimal homeomorphism $\f$ of a Cantor set
$X$ and involutions from $[\f]$ with clopen supports  belong to $D([\f])$.
This will allow us to prove that the simple group $D([\f])$ is a complete
invariant for the class of orbit equivalence of $\f$. As a corollary of
the technique used, we also get that $\f$ can be written as a product of three
involutions from $[\f]$. Our considerations are mainly based on the ideas
of Fathi \cite{fathi}.

\begin{remark}\label{PeriodicCommutator} Suppose that $g$ is a periodic
homeomorphism from $[\f]$ such that the space $X$ can be decomposed into
clopen sets
$$X=\bigcup_{i\in I}\bigcup_{j=0}^{n_i-1}g^j(X_i^0)$$
with $g^{n_i}(x)=x$ for all $x\in X_i^0$.  We give two model
situations when $g$ can be easily written as a commutator in
$[\f]$. Consider the following cases:

 (1) Suppose that $n_i$ is odd for $i\in I$. Then $g|X_i$ is a
 commutator in $[\f]$. Indeed,
 let $m=\frac{n_i-1}{2}$ and define the homeomorphisms $g_1$ and $g_2$ as
 follows: $$g_1(x)=\left\{\begin{array}{lll}g(x) &\mbox{if }x\in
 \bigcup_{k=0}^{m-1}g^k(X_i^0)\\
                 g^{-m}(x)&\mbox{if }x\in g^m(X_i^0)\\
                 id &\mbox{elsewhere}.     \end{array}\right.$$

 $$g_2(x)=\left\{\begin{array}{lll}g(x) &\mbox{if }x\in
 \bigcup_{k=m}^{n_i-2}g^k(X_i^0)\\
                 g^{-m}(x)&\mbox{if }x\in g^{n_i-1}(X_i^0)\\
                 id&\mbox{elsewhere}.     \end{array}\right.$$
Then $g=g_1g_2$.  Since $g_2=\psi g_1^{-1}\psi^{-1}$ for some $\psi\in [\f]$,
we see that $g=g_1g_2=[g_1,\psi]$.

(2) Now suppose that $n_i$ is even and $X_i^0=X_i^0(l)\sqcup X_i^0(r)$ where
$X_i^0(l)$ and $X_i^0(r)$ are $[\f]$-equivalent clopen sets, i.e. $\alpha(X_i^0(l))
=X_i^0(r)$ for some  $\alpha\in [\f]$.
Then $g|X_i$ can be written down as a commutator in  $[\f]$. Indeed, put $X_i(l)=\bigcup_{k=0}^{n_i-1} g^k (  X_i^0(l))$ and $X_i(r)=\bigcup_{k=0}^{n_i-1}
g^k( X_i^0(r))$. Define   $l,r\in [\f]$ as follows: $l|X_i(l)=g|X_i(l)$ and $l=id$
elsewhere; $r|X_i(r)=g|X_i(l)$ and $r=id$  elsewhere. Note that $g=lr$.
By Theorem \ref{equivalent-clopen-sets}, it is easy to see that
$r=\psi l^{-1}\psi^{-1}$ for some  $\psi\in [\f]$. Hence   $g=lr=[l,\psi]$.
\end{remark}

The proof of the fact that $\f\in D([\f])$ consists of series of lemmas.

\begin{lemma}\label{MinimalFirstStep} Suppose that $f\in [\f]$ has clopen support
and $f|supp(f)$ is minimal. Then for given $\delta > 0$ there exist
$f_1,s,t\in [\f]$ such that

(1) $f=f_1[s,t]$;

(2) $supp(f_1)$ is clopen, $f_1|supp(f_1)$ is minimal, and
$\mu(supp(f_1))<\delta$ for all $\mu\in \M$;

(3) $supp(s)\cup supp(t)\cup supp(f_1)\subseteq supp(f)$.
\end{lemma}
{\it Proof.}  Take any clopen set $ A\subset supp(f)$ such that
$\mu(A)<\delta/2$ for all $\mu\in \M$ and $f^i(A)\cap A=\emptyset$ for
$i=1,2,3$.

Applying the first return function, construct the K-R partition $\Xi$ over
$A$ (see Section \ref{Section-Preliminaries} for the definition). Suppose
that $\Xi=\{\xi_1,\ldots,\xi_n,\xi_1',\ldots,\xi_m'\}$ where $\xi_i$ are
the $f$-towers  with even heights and  $\xi_j'$ are $f$-towers with odd
heights. Let $h(\xi)$ denote the height and $B(\xi)$ denote the base of an
$f$-tower $\xi$. Set
$$B=A\cup\bigcup_{i=1}^nf^{h(\xi_i)/2}B(\xi_i).$$
Clearly, $\mu(B)<\delta$  for all $\mu\in \M$.

Define $f_1$ as the induced homeomorphism $f_B$. Again using the first
return function, consider the K-R-partition $\mathcal P$ over $B$.
Note that $\mathcal P=\{\xi_1',\ldots,\xi_m',\ \xi_1^1,\
\xi_1^2,\ldots,\xi_n^1,\\ \xi_n^2\}$, where the $\xi_i'$'s are $f$-towers
with odd heights as above, $\xi_i^1$ is the lower half of $\xi_i$, and
$\xi^2_i$ is the upper half.

Define the periodic homeomorphism $g$ as follows:
$$g(x)=\left\{\begin{array}{ll}f(x)&\mbox{if } x\notin \bigcup_{\xi\in\mathcal P}f^{h(\xi)-1}(B(\xi)),\\
\\
  f^{-h(\xi)+1}(x) & \mbox{if }x\in f^{h(\xi)-1}(B(\xi))\mbox { for some }\xi\in\mathcal P.\end{array}\right.$$
 Then $f=f_1g$. Since $\xi_i'$ has odd height,  Remark  \ref{PeriodicCommutator} implies that
 $g| \xi_i'$  is a commutator in $[\f]$.
Consider $g|(\xi_i^1\cup\xi_i^2),\ i=1,...,n$. By  construction, the bases
$B(\xi_i^1)$ and $B(\xi_i^2)$ are $[\f]$-equivalent. Therefore, the
application of Remark \ref{PeriodicCommutator}(2) ensures us that
$g=[s,t]$ for some $s,t\in [\f]$. Statements (2) and (3) are obvious.
\hfill$\square$

%
%
\begin{remark}\label{notat}  Clearly, we can construct $f_1$ in Lemma
\ref{MinimalFirstStep} such that $supp(f_1)$ is always a proper subset of
$supp(f)$. Let now $x_0$ be a point from $supp(f)\setminus supp(f_1)$.
Then one can find a sequence $\{C_n\}_{n\geq 1}$ of mutually disjoint
clopen sets such that (1) $C_i\subseteq supp(f)$ (2) $C_1=supp(f_1)$; (3)
$x_0\notin C_n$ for $n\geq 1$; (4) $diam(C_n\cup\{x_0\})\to 0$ as
$n\to\infty$. Put $\delta_n=\inf\{\mu(C_n) : \mu\in \M\}$. By Proposition
\ref{propertiesOfInvariantMeasures}, every $\delta_n>0$.
\end{remark}

\begin{lemma}\label{DistributionOfMinimal} Suppose that $f_1$ and a sequence
$\{C_n\}$ are as in Remark \ref{notat}. Then there exists two sequences
$\{f_i\}_{i\geq 1}$ and $\{g_i\}_{i\geq 1}$ of homeomorphisms from $[\f]$
such that for every $i\geq 1$

(1)  $supp(f_i)$ is clopen and $f_i$ is minimal on it;

(2) $f_i=f_{i+1} g_i$;

(3) $supp(f_i)\subseteq C_i$;

(4) $g_i=[s_i',t_i'][s_i,t_i]$ for some $s_i,t_i,s_i',t_i'\in [\f]$ with
$supp(s_i)\cup supp(t_i)\cup supp(s_i')\cup supp(t_i')\subseteq C_i\cup
C_{i+1}$.
\end{lemma}

{\it Proof.} We will explain only the first step of the construction:
given $f_1, C_1$, and $C_2$, we will find  $f_2, g_1$.

By Lemma \ref{MinimalFirstStep}, there exist homeomorphisms $\overline
f,s,t\in [\f]$ with clopen supports in the set $C_1$ such that $f_1=
\overline f [s,t]$ and $\mu(supp(\overline f))<\delta_2$ for all $\mu\in
\M$. It follows from Theorem \ref{equivalent-clopen-sets} that there
exists a homeomorphism $t'\in [\f]$ with clopen support in the set
$C_1\cup C_2$ such that $t'(supp(\overline f))\subset C_2$. Set
$f_2=t'\overline f t'^{-1}$. Then the support of $f_2$ is a clopen subset
of $C_2$ and  $f_2\in [\f]$ is minimal on $supp(f_2)$. Thus, we have that
$$f=\overline f [s,t]= f_2 [t',\overline f^{-1}][s,t].
$$
Setting $g_1=[t',\overline f^{-1}][s,t]$, we complete the proof.
\hfill$\square$

\begin{lemma}\label{perfectnessMinimal}  Let $f_1$ be as in Lemma
\ref{MinimalFirstStep}. Then  $f_1\in D([\f])$ and $f_1$ is a product of
four commutators from $[\f]$.
\end{lemma}
{\it Proof.} Let  $\{f_i\}_{i\geq 1}$ and $\{g_i\}_{i\geq 1}$ be the
sequences of homeomorphisms constructed in Lemma
\ref{DistributionOfMinimal}. Recall that $g_i=[s'_i,t'_i][s_i,t_i]$ and
$supp(s_i')\cup supp(t_i')\cup supp(s_i)\cup supp(t_i)\subseteq C_i\cup
C_{i+1}$. Note that the homeomorphisms $\{g_{2k+1}\}_{k\geq 0}$ have mutually
disjoint  supports. So do $\{g_{2k}\}_{k\geq 1}$. Define the maps
$g_{odd}$ and $g_{even}$ as follows:
$$
g_{odd}(x)=\left\{\begin{array}{ll} g_i(x) & \mbox{whenever }x\in supp(g_i)\mbox{ for odd }i'\\
x & \mbox{elsewhere} \end{array}\right.
$$
$$
g_{even}(x)=\left\{\begin{array}{ll} g_i(x) & \mbox{whenever }x\in supp(g_i)\mbox{ for even }i\\
x & \mbox{elsewhere.} \end{array}\right.
$$
It follows from the choice of the sets $C_i$ and the property
$diam(C_n\cup\{x_0\})\to 0$ that $g_{odd}$ and $g_{even}$ are
homeomorphisms.

We see $g_{odd}=[s_{odd}',t_{odd}'][s_{odd},t_{odd}]$ and
$g_{even}=[s_{even}',t_{even}'][s_{even},t_{even}]$ where the homeomorphisms
$s_{odd}'$, $s_{even}'$, $t_{odd}'$, $t_{even}'$ $s_{odd}$, $s_{even}$,
and $t_{odd}$, $t_{even}$ are defined similarly to $g_{odd}$ and
$g_{even}$.

By definition of $g_i$, we have that $g_i=f_{i+1}^{-1} f_i $. Since all the
$f_i$'s have disjoint supports, we can formally write down the infinite products
$$g_{odd}=(f_2^{-1} f_1) (f_4^{-1} f_3) (f_6^{-1} f_5)\cdots= (f_1f_3f_5\cdots)(f_2^{-1}f_4^{-1}f_6^{-1}\cdots),$$
$$g_{even}=(f_3^{-1}f_2)(f_5^{-1}f_4)(f_7^{-1}f_6)\cdots=
(f_3^{-1}f_5^{-1}f_7^{-1}\cdots)(f_2f_4f_6\cdots).$$
Therefore, $f_1 =g_{odd}g_{even}=[s_{odd}',t_{odd}'][s_{odd},t_{odd}]
[s_{even}',t_{even}'] [s_{even},t_{even}]$. \hfill$\square$

\begin{theorem}\label{ProductCommutators} Let $(X,\f)$ be a Cantor minimal
system. Suppose that a homeomorphism $f\in [\f]$ has clopen support and
is minimal on it. Then $f\in D([\f])$ and $f$ is a product of  five commutators
from $[\f]$. In particular, $\f\in D([\f])$.
\end{theorem}
{\it Proof.} The proof follows immediately from Lemmas
\ref{MinimalFirstStep} and \ref{perfectnessMinimal}. \hfill$\square$

\begin{remark}\label{TopolSimplicity}  Let $\mathcal P_1(X)$ denote the
set of all Borel probability measures on $X$. For any $g\in Homeo(X)$,
$\e>0$, and any $\mu_1,\ldots,\mu_n\in \mathcal P_1(X)$, define
$U(g;\mu_1,\ldots,\mu_n;\e)=\{h\in Homeo(X) : \mu_i(\{x\in X : h(x)\neq
g(x)\})<\e\mbox{ for }i=1,\ldots,n\}$. Let $\tau$ denote the topology on
$Homeo(X)$ generated by the base sets of the form
$U(g;\mu_1,\ldots,\mu_n;\e)$. This topology was defined and studied in
\cite{bezuglyi_dooley_kwiatkowski}. In particular, it was shown that the
topological group $[\f]$ is closed in $Homeo(X)$
with respect to $\tau$. On the other
hand the topological full group $[[\f]]$ is $\tau$-dense in $[\f]$
\cite{bezuglyi_kwiatkowski} (see also \cite{medynets:2005} for another
proof of the result).

Developing the ideas used in this section, one can show that
$[[\f]]\subset D([\f])$.  Hence, Theorem \ref{simplicity} implies that the
full group of a Cantor minimal system has no  $\tau$-closed normal
subgroups.
\end{remark}

Now, we present several immediate corollaries of Theorem
 \ref{ProductCommutators}.
\begin{corollary}\label{ClopenInvol} Suppose $\f$ is a minimal homeomorphism
and  $g\in [\f]$ is an involution with
clopen support. Then $g$ is a product of $10$ commutators in $[\f]$.
\end{corollary}
{\it Proof.} By Lemma \ref{PeriodicProperties}, there exists a clopen set
$A\subset supp(g)$ such that $supp(g)$ is a disjoint union of
$A$ and $g(A)$. Define $f_1$ as the induced map $\f_A$ and let $f_2=f_1g$.
Clearly, $f_1$ and $f_2$ have clopen supports and are minimal on their
supports. Then Theorem \ref{ProductCommutators} asserts that
$g=f_1^{-1}f_2$ is a product of 10 commutators from $[\f]$.
\hfill$\square$

\smallskip
It follows from Theorems \ref{ProductCommutators} and \ref{simplicity}
that $\f\in H$ for any normal subgroup $H$ of $[\f]$. In particular, this
fact allows us to show that $\f$ is a product of involutions. The
following result gives the number of involutions needed to represent $\f$.
The proof is based on Theorem \ref{simplicity}.

\begin{corollary}\label{Prod18Invol} A minimal homeomorphism
$\f$ is a product of 18 involutions from
$[\f]$ which have clopen supports.
\end{corollary}
{\it Proof.} Let $w$ be any involution from $[\f]$ with clopen support.
Choose any clopen set $E$ such that $w(E)\cap E=\emptyset$. Set
$\delta=\inf\{\mu(E) : \mu\in\M\}>0$. Take any clopen set $A'$ with
$\mu(A')<\delta$ for all $\mu\in \M$. Let $\f_A$ be the induced
homeomorphism of $\f$, where $A$ is a proper clopen subset of $A'$.  It follows from Proposition \ref{InducedMap} that
$g=\f_A^{-1}\f$ is a periodic homeomorphism from $[[\f]]$. It is not hard
to see that $g=st$, where $s$ and $t$ are involutions from $[[\f]]$ (see,
for example, \cite[Proposition 4.1]{miller:thesis}). Applying Lemma
\ref{perfectnessMinimal} to $\f_A (=f_1)$, we obtain  that
$\f_A=[s_1,t_1]\ldots [s_4,t_4]$, where $supp(s_i)\cup supp(t_i)\subseteq
A'$ and $s_i,t_i\in [\f]$ for each $i=1,\ldots,4$.

We claim that if some homeomorphisms $h$ and $g$ are supported on $A'$,
then $[h,g]$ is a product of four conjugates of $w$. Indeed, by Theorem
\ref{equivalent-clopen-sets}, find $\alpha\in [\f]$ such that
$\alpha(A')\subseteq E$. Denote $q =\alpha^{-1}w\alpha$. Clearly, $q(A')\cap
A'=\emptyset$. Observe that $qhq^{-1}$ commutes with $h$ and $g$. Hence
$$\begin{array}{llll}[h,g]& = & hgh^{-1}g^{-1}\\ &= & h(qh^{-1}q^{-1})(qhq^{-1})
gh^{-1}g^{-1}\\ &= & h(qh^{-1}q^{-1})g(qhq^{-1})h^{-1}g^{-1}\\& = &
(hqh^{-1}) (q^{-1})(gq g^{-1})(gh q^{-1}h^{-1}g^{-1}).
\end{array}$$
This implies that $[h,g]$ is the product of four conjugates of $w$.
Therefore, $\f_A$ is a product of 16 conjugates of $w$ and $\f= \f_Ag$ is
a product of 18 involutions. \hfill$\square$

\medskip
We can refine Corollary \ref{Prod18Invol} and show that $\f$ can be
written, in fact, as a product of three involutions with clopen supports
from $[\f]$. In our proof we follow the arguments from
\cite{ryzjikov:1985}. Recall that everywhere below $\f$ denotes a minimal
homeomorphism.
%
%

\begin{definition} We say that a homeomorphism $g\in [\f]$ is an
{\it $n$-cycle} on
disjoint clopen sets $E_0,E_1,\ldots,E_{n-1}$  if: (1) $g(E_i)=E_{i+1}$ for
$i=0,\ldots,n-2$ and $g(E_{n-1})=E_0$; (2) $g(x)=x$ for all $x\in
X\setminus(E_0\cup\ldots\cup E_{n-1})$.
\end{definition}

\begin{lemma}\label{CyclesLemma} Let $g\in [\f]$ be an $18$-cycle on disjoint
clopen sets $E_0,E_1,\ldots, E_{17}$ such that $g^{18}|E_0$ is a minimal
homeomorphism. Then there exists an involution $d\in [\f]$ with clopen
support such that the homeomorphism $gd$ is an $18$-cycle on $E_0,\ldots,
E_{17}$, and $(gd)^{18}=id$.
\end{lemma}
{\it Proof.} By  Corollary \ref{Prod18Invol}, there are involutions
$h_0,\ldots,h_{17}$ from $[\f]$ with clopen supports such that
$g^{18}|E_0=h_0\ldots h_{17}$  and $supp(h_i)\subseteq E_0$.

Set $d_k=g^{k}h_k^{-1}g^{-k}$ and $d=d_0d_1\ldots d_{17}$. Then $d$ is an
involution with clopen support and
$$\begin{array}{llll}(gd)^{18}|E_0 &=&
 (gd_{17})\ldots(gd_1)(gd_0)|E_0\\ &= &(g g^{17}h_{17}^{-1}g^{-17})(g
g^{16}h_{16}^{-1}g^{17})\ldots (g g^2 h_2^{-1}g^{-2})
 (g gh_1^{-1}g^{-1}) (g h_0^{-1})|E_0 \\&=&
 g^{18}h_{17}^{-1}\ldots h_0^{-1}|E_0\\ &=&id.
 \end{array}$$
 Since $(gd)$ is an $18$-cycle on $E_0,\ldots,E_{17}$, we get that
 $(gd)^{18}=id$. \hfill$\square$

\begin{remark}\label{InvolutionRemark} Since the homeomorphism $gd$ has the
period $18$ on its support, i.e. for all $x\in supp(gd)$ one has
$(gd)^i(x)=x$ iff $i=18k$ for  $k\in\mathbb Z$,  there are two involutions
$s,t\in [\f]$ such that $gd=st$.   Furthermore, $supp(s) \subseteq E_1\cup
E_2\cup\ldots E_{17}$.
\end{remark}

 \begin{lemma}\label{TowerLemma} Given $n> 1$, there exists
 a clopen set $A$ such that

 (1) $A,\f (A),\ldots, \f^{n-1}(A)$ are disjoint;

 (2) $\f (B)\subseteq A$, where $B=X\setminus \bigcup_{i=0}^{n-1}\f^i(A)$.

 \end{lemma}
{\it Proof.} Let $E$ be a clopen set such that $\f^i(E)\cap E=\emptyset$
for $i=1,\ldots,n^2$. Construct a K-R partition over $E$, i.e.
$E=\bigcup_kE_k$, where $E_k=\{x\in E : \f^k(x)\in E \mbox{ and
}\f^j(x)\notin E\mbox{ for }0<j<k\}$.  Observe  that for all $k\leq n^2$,
we have that $E_k=\emptyset$.

(1) If $k=nl+r$, $0< r < n$, then we set
$$
A_k=\bigcup_{i=0}^{r-1}\f^{(n+1)i}(E_k)\cup\bigcup_{i=r}^{l-1}\f^{ni+r-1}(E_k).
$$
This means that we choose $r$ times every $(n+1)$-th set from the family
$\{E_k,\f(E_k),\ldots,\\ \f^{k-1}(E_k)\}$ starting from the first one and
then we take every $n$-th set.

(2) If $k=nl$, then we set
$$A_k=\bigcup_{i=0}^{l-1}\f^{ni}(E_k).$$

Denoting $A=\bigcup_k A_k$, we get the result. \hfill$\square$

\begin{theorem}\label{Theorem_Product3involutions} Let $(X,\f)$ be a
Cantor minimal system. Then there exist involutions $i_1,i_2,i_3\in [\f]$
with clopen supports  such that $\f=i_1i_2i_3$.
\end{theorem}
{\it Proof.} Find a clopen set $A$ satisfying the conditions of Lemma
\ref{TowerLemma} for $n=18$. Let $B=X\setminus \bigcup_{i=0}^{17}\f^i(A)$.
Clearly, $\f(B)\cap B=\emptyset$. Define an involution $b$ as follows
$b|B=\f|B$ and $b|\f(B)=\f^{-1}|\f(B)$. It follows that $g=b\f$ is an
$18$-cycle on the clopen sets $A_0,\ldots, A_{17}$ where $A_i=\f^{i}(A)$.

By Lemma \ref{CyclesLemma}, we can find an involution $d$ with clopen
support such that $(gd)^{18}=id$. It follows from Remark
\ref{InvolutionRemark} that there exist involutions $s$ and $t$ such that
$gd=st$. This implies that $\f=b^{-1}st d^{-1}$ is the product of four
involutions. As mentioned in Remark \ref{InvolutionRemark},
$supp(s)\subseteq A_1\cup\ldots A_{17}$. Hence $supp(s)\cap
supp(b)=\emptyset$. It follows that $w=b^{-1}s$ is an involution. This
proves that $\f=wtd^{-1}$ is the product of three involutions with clopen
supports. \hfill$\square$

%
%

\section{Flip conjugacy and orbit equivalence}\label{SectionFlipConjugacy}
 In the section, we show that the classes of orbit equivalence and flip
 conjugacy of a Cantor minimal system are completely determined by simple groups.

 \begin{definition}   (1) Cantor minimal
 systems $(X_1,\f_1)$ and $(X_2,\f_2)$    are called {\it orbit equivalent}
 if there exists a homeomorphism $F: X_1\rightarrow X_2$ such that
 $F(Orb_{\f_1}(x))=Orb_{\f_2}(F(x))$ for all $x\in X_1$.

 (2) $(X_1,\f_1)$ and $(X_2,\f_2)$ are called {\it flip conjugated} if there exists a
 homeomorphism $F: X_1\rightarrow X_2$ such that $F\circ\f_1\circ F^{-1}$
 is equal to either $\f_2$ or $\f_2^{-1}$.
 \end{definition}

 The following theorem can be deduced from  \cite{gps:1999}.

 \begin{theorem}\label{Theorem_OrbitEquivalence} Let $(X_1,\f_1)$ and $(X_2,\f_2)$
 be Cantor minimal systems. The homeomorphisms $\f_1$ and $\f_2$ are orbit
 equivalent if and  only if $D([\f_1])$ and $D([\f_2])$ are isomorphic as algebraic
 groups.
 \end{theorem}
{\it Proof.} The theorem can be proved in the same way as Corollary 4.6
from \cite{gps:1999}. We notice only that according to Corollary
\ref{ClopenInvol} every involution with clopen support belongs to the
commutator subgroup. We recall also that $\f\in D([\f])$ (see Theorem
\ref{ProductCommutators}). \hfill$\square$

\medskip Now, we start studying the class of flip conjugacy of a Cantor minimal
system in terms of topological full groups. In our arguments, we mainly
follow the proof of Theorem 384D from \cite{fremlin}. This theorem was
used to show that certain transformation groups of complete Boolean
algebras have no outer automorphisms (see also \cite{eigen:1982}).

 Let $\Gamma$  denote one of
the following groups: (1) the topological full group $[[\f]]$;  (2) the
group $[[\f]]_0$ defined in Remark \ref{Remark_NormalSubgroups}; (3) the
commutator subgroup $D([[\f]])$ of $[[\f]]$. Notice that we have the
following inclusion:
\begin{equation}\label{threeSubgroups}
D([[\f]])\subseteq [[\f]]_0\subsetneq [[\f]].
\end{equation}

We will show that every group from (\ref{threeSubgroups}) is a complete
invariant of flip conjugacy of $\f$. The proposed proof  works for any of
these groups $\Gamma$, because the only group property we exploit is the
existence  of many `small involutions' in the following sense.

\begin{lemma}\label{ManyInvolutions} Each group $\Gamma$ has many
involutions, in the sense that for any clopen set $A$ and any $x\in A$,
there exists $h\in \Gamma$ such that $hx\neq x$, $h^2=1$ and
$supp(h)\subseteq A$. Moreover, for every $n> 0$, there exists $h\in
\Gamma$ such that $h$ is supported by $A$, $x\in supp(h)$ and $h|supp(h)$
has period $n$.
\end{lemma}
{\it Proof.} By (\ref{threeSubgroups}), it suffices  to establish the
result for $D([[\f]])$ only. Let us find integers
$0=m_0<m_1<\ldots<m_{2n-1}$ such that $\f^{m_i}(x)\in A $ for $i=0,\ldots
2n-1$. Take a clopen neighborhood $V$ of $x$ such that $\f^{m_i}(V)\subset
A$  and $\f^{m_i}(V)\cap \f^{m_j}(V)=\emptyset$ for all $i,j=0,\ldots,
2n-1$, $i\neq j$.

Define a homeomorphism $l\in [[\f]]$ as follows:
$l(x)=\f^{m_{i+1}-m_i}(x)$ if $x\in \f^{m_i}(V)$ for $i=0,1,\ldots,
n-2$, $l(x)=\f^{-m_{n-1}}(x)$ if $x\in \f^{m_{n-1}}(V)$, and $l(x)=x$
elsewhere. Similarly, define a homeomorphism $r\in [[\f]]$:
$r(x)=\f^{m_{i+1}-m_i}(x)$ if $x\in \f^{m_i}(V)$ for $i=n,n+1,\ldots,
2n-2$,  $r(x)=\f^{-(m_{2n-1}-m_n)}(x)$ if $x\in \f^{m_{2n-1}}(V)$, and
$r(x)=x$ elsewhere.

 It is not hard to see that there
exists $\alpha\in [[\f]]$ such that $l=\alpha r\alpha^{-1}$.
Therefore, $h=lr^{-1}\in D([[\f]])$ and $h$ has period $n$ on its support $\f^{m_0}(V)\cup\ldots \cup \f^{m_{2n-1}}(V)$.
\hfill$\square$

\medskip As a corollary of the lemma, we obtain the following result.

\begin{corollary}\label{TopologyGeneration} The family of clopen sets
$\{supp(g) : g\in \Gamma\mbox{ and }h^2=1\}$ generates the clopen topology
of $X$.
\end{corollary}

We will need some notions of the theory of Boolean algebras. We refer the
reader to the book \cite{fremlin} for a comprehensive coverage of the
theory of Boolean algebras and their automorphisms.

Let $X$ be a Cantor set. Recall that an open set $A$ is called {\it
regular open} if $A=int(\overline{A})$. Denote the family of all regular
open sets by $RO(X)$. Notice that the family of clopen sets $CO(X)$ is
contained in $RO(X)$.

Let $\mathcal A$  be a Boolean algebra and $\mathcal H\subset \mathcal A$.
Define $sup(\mathcal H)$ to be the smallest element of $\mathcal A$ that
contains all elements of $\mathcal H$. If $sup(\mathcal H)$ exists for any
family $\mathcal H\subseteq \mathcal A$, then the Boolean algebra
$\mathcal A$ is called {\it complete}.

\begin{theorem}[Theorem 314P of \cite{fremlin}]
$RO(X)$ is a complete Boolean algebra with Boolean operations given by
$$A\bigvee B= int (\overline{A\cup B}),\;\; A\bigwedge B= A\cap B,\;\;
A\setminus_{RO(X)} B=A\setminus \overline{B}$$
and with suprema given by $sup(\mathcal H)=int(\overline{\bigcup \mathcal
H })$.
\end{theorem}

\begin{remark} Notice that  finite theoretical-set unions of clopen sets coincide
with the Boolean ones.
\end{remark}

\begin{lemma}\label{ClopenSet} If $A,B\in RO(X)$ and $A\not\subset B$, then $A\setminus
B$ contains a clopen set.
\end{lemma}

{\it Proof.} Let $C=A\cap B$. If $C=\emptyset$, then the result is clear.
Assume that  $C\neq\emptyset$ and $A\setminus C=A\setminus B$ has no
internal points. Then $A\subset \overline{C}$. This implies that $A=C$, a
contradiction. \hfill$\square$

\medskip
Now we are ready to start the proof of the main result of this section.
Recall that for a Cantor minimal system $(X,\f)$,  $\Gamma$ stands for one
of the groups: (1) $[[\f]]$, (2) $[[\f]]_0$,  (3) $D([[\f]])$.

%
%

\begin{theorem}\label{Theorem_Spatial_Realization} Let $(X_1,\f_1)$ and
$(X_2,\f_2)$ be Cantor minimal systems. If $\alpha: \Gamma_1\rightarrow \Gamma_2$
is a group isomorphism, then $\alpha$  is spatial, i.e. there exists a
homeomorphism $\widehat \alpha : X_1\rightarrow X_2$ such that $\alpha(g)
=\widehat\alpha g\widehat\alpha^{-1}$ for any $g\in \Gamma_1$.
\end{theorem}

{\it Proof.} The idea of the proof is the following: we study the local
subgroup $\Gamma_A := \{g\in \Gamma : gx=x \mbox{ for all } x\in X\setminus A \}$,
where a clopen set $A$  is the support of an
involution from $\Gamma$, and  describe $\Gamma_A$ in group terms. This
description allows us to construct an automorphism $\widehat{\alpha} :
RO(X_1)\rightarrow RO(X_2)$ of Boolean algebras of regular open sets,
which also sends clopen sets onto clopen sets. This automorphism gives
rise to the spatial realization of $\alpha$. For convenience, we will omit
the index $i$ in the notation of Cantor minimal system $(X_i,\f_i)$ and
the group $\Gamma_i$. We split the proof of the theorem into two lemmas.

\begin{definition} Let $\pi\in \Gamma$ be any involution. Set
$$
C_\pi=\{g\in \Gamma : g\pi=\pi g \},
$$
the centralizer of $\pi$ in $\Gamma$;
$$
U_\pi=\{g\in C_\pi : g^2=1\mbox{ and }g(hgh^{-1})=
(hgh^{-1})g\mbox{ for all }h\in C_\pi\},
$$
the involutions from $C_\pi$ which commute with all their conjugates in
$C_\pi$;
$$
V_\pi=\{g\in \Gamma : gh=hg\mbox{ for all }h\in U_\pi\},
$$
the centralizer of $U_\pi$ in $\Gamma$;
$$S_\pi=\{g^2 : g\in V_\pi\};
$$
and
$$
W_\pi=\{g\in \Gamma  : gh=hg\mbox{ for all }h\in S_\pi\},
$$
the centralizer of $S_\pi$ in $\Gamma$.
\end{definition}
\medskip

Clearly, for any involution $\pi$, $supp(\pi)$ is a clopen set and
$\alpha(W_\pi)=W_{\alpha(\pi)}$ where $\alpha$ is as in Theorem \ref{Theorem_Spatial_Realization}.

\begin{lemma}\label{Fremlin_Lemma1} Let $\pi\in \Gamma$ be an involution. Then
$W_\pi=\Gamma_{supp(\pi)}$.
\end{lemma}
{\it Proof.} To prove the result, we will consequently study
 the properties of $C_\pi$, $U_\pi$, $V_\pi$, $S_\pi$, and  $W_\pi$.

\medskip (1) $g(supp(\pi))=supp(\pi)$ for all $g\in C_\pi$. $\vartriangleleft$
It is  easy to see that $supp(g\pi g^{-1}) = g(supp(\pi))$. Since $g\pi
g^{-1} = \pi$, one has that $g(supp(\pi))=supp(\pi)$. $\vartriangleright$

\medskip (2-i)  $supp(g)\subseteq supp(\pi)$ for all $g\in U_\pi$.
$\vartriangleleft$ Assume the converse. Then there exists a clopen set
$A\subset X\setminus supp(\pi)$ such that $gA\cap A=\emptyset$. By Lemma
\ref{ManyInvolutions}, find a homeomorphism $h\in \Gamma$ with support in
$A$ such that for a clopen set $V\subset A$ one has that $h^i(V)\cap
V=\emptyset$, $i=1,2$. Note that $h\in C_\pi$. Then
$$g(hgh^{-1})(V)=g^{2}h^{-1}(V)=h^{-1}(V)\mbox{ whereas }(hgh^{-1})g(V)=hg^2(V)=h(V).$$
The choice of $h$ guarantees us  that $g(hgh^{-1})\neq (hgh^{-1})g$. Hence
$g\notin U_\pi$, which is a contradiction. $\vartriangleright$

(2-ii) If a clopen set $A$ is $\pi$-invariant, then $\pi_A\in U_\pi$
where the homeomorphism $\pi_A$ coincides with $\pi$ on $A$ and is equal
to $id$ elsewhere.

\medskip (3-i) $V_\pi\subset C_\pi$, because $\pi\in U_\pi$.

(3-ii) If $g\in V_\pi$, then $g(B)\subseteq B\cup \pi(B)$ for all clopen
sets $B\subseteq supp(\pi)$. $\vartriangleleft$ Assume the converse, i.e.
$g(B)$ is not in $B_0=B\cup \pi(B)$ for some set $B$. Note that
$\pi(B_0)=B_0$. As $B\subset B_0$, we have that $C=g(B_0)\setminus
B_0\neq\emptyset $. Since $\pi g(B_0)=g\pi(B_0)=g(B_0)$, we obtain that
$$\pi(C)=\pi(g(B_0)\setminus B_0)=\pi g(B_0)\setminus
\pi(B_0)=g(B_0)\setminus B_0=C.
$$
Since $g\in V_\pi\subset C_\pi$, we have that $supp(\pi)$ is $g$-invariant
and  $C\subset supp(\pi)$. As $\pi$ is an involution, then $C=C'\sqcup
C''$ for some clopen set $C'$ with $C''=\pi(C')$. Note that $g(C)\cap
C=\emptyset$. Therefore
$$\pi_Cg(C')=g(C')\neq g(C'') =g\pi_C(C').
$$
Thus, we obtain that $g$ does not commute with $\pi_C\in U_\pi$,  a
contradiction. $\vartriangleright$

(3-iii) If $g\in V_\pi$, then $g^2(B)=B$ for any clopen $B\subseteq
supp(\pi)$. $\vartriangleleft$ Assume the converse, i.e. there exists a
clopen set $B\subset supp(\pi)$ such that $g^2(B)\cap B=\emptyset$. We can
also  assume that $g(B)\cap B=\emptyset$. We know that  $g(B)\subset B\cup
\pi(B)$. This implies that $g(B)\subset \pi(B)$. By the same argument
applied to $g(B)$, we obtain that $g^2(B)\subset \pi^2(B)=B$. If
$g^2(B)\neq B$, then $\mu(B)=0$ for any $\mu\in \M$, which contradicts to
the minimality of $\f$. Therefore, $g^2(B)=B$. $\vartriangleright$

\medskip (4-i) If $g\in S_\pi$, then $supp(g)\subset X\setminus supp(\pi)$.
$\vartriangleleft$ It follows  from  (3-iii). $\vartriangleright$

(4-ii) For any clopen set $C\subset X\setminus supp(\pi)$, there is an
involution $h\in S_\pi$ supported on $C$. $\vartriangleleft$ By Lemma
\ref{ManyInvolutions}, there exists a periodic homeomorphism $g$ of order
$4$ with support in $C$. Property (2-i) implies that $g\in V_\pi$.
Therefore, $g^2\in S_\pi$. $\vartriangleright$

\medskip (5) $W_\pi=\Gamma_{supp(\pi)}$. $\vartriangleleft$  It follows from
 (4-i)  that $\Gamma_{supp(\pi)}\subset W_\pi$. To get the converse inclusion, we
 consider any $g\in W_\pi$ and suppose that there exists a clopen set
 $B\subset X\setminus supp(\pi)$
such that $g(B)\cap B=\emptyset$. By  (4-ii), find an involution $h$ from
$S_\pi$ with support in $B$. Take any clopen set $C\subset B$ with
$h(C)\cap C=\emptyset$. Therefore,
$$hg(C)=g(C)\neq gh(C).$$
This implies that $gh\neq hg$, a contradiction. Thus,
$W_\pi=\Gamma_{supp(\pi)}$. $\vartriangleright$

This completes the proof of the lemma. \hfill$\square$

\medskip The following lemma gives the spatial realization of the
group isomorphism $\alpha$.

\begin{lemma} Let $(X_i,\f_i), \Gamma_i, \alpha$ be as in Theorem
\ref{Theorem_Spatial_Realization}. The map $\Lambda: RO(X_1)\rightarrow
RO(X_2)$ given by
$$\Lambda(A)=\bigvee \{supp(\alpha(\pi)) : \pi\in \Gamma_1,\;\pi^2=
1\mbox{ and }supp(\pi)\subseteq A\}$$
is a Boolean algebra isomorphism. Furthermore, $\Lambda
(CO(X_1))=CO(X_2)$.
\end{lemma}

{\it Proof.} (1) Note that if $g,\pi\in \Gamma_1$ and $\pi^2=1$, then
$$supp(g)\subseteq supp(\pi) \Leftrightarrow supp(\alpha(g))
\subseteq supp(\alpha(\pi)).
$$
$\vartriangleleft$ By Lemma \ref{Fremlin_Lemma1} (5), we see that
$supp(g)\subseteq supp(\pi)$ iff $g\in W_\pi$ iff $\alpha(g)\in
W_{\alpha(\pi)}$ iff $supp(\alpha(g)) \subseteq supp(\alpha(\pi))$.
$\vartriangleright$

\medskip (2) It easily follows from the definition of $\Lambda$ that
$\Lambda$ is order-preserving, i.e. if $A\subseteq B$
with $A,B\in RO(X_1)$, then $\Lambda(A)\subseteq \Lambda(B)$.

\medskip (3) Let $\pi\in \Gamma_1$ be such that $\pi^2=1$ and $A\in RO(X_1)$. If
$supp(\pi)\not\subset A$, then $supp(\alpha(\pi))\not\subset \Lambda(A)$.
$\vartriangleleft$  Take a non-empty clopen set $V\subseteq
supp(\pi)\setminus A$ and a homeomorphism $h\in \Gamma_1$ of order $4$
with support in $V$ (see Lemma \ref{ClopenSet}). As $supp(h)\subset
supp(\pi)$, we have $supp(\alpha(h))\subset supp(\alpha(\pi))$ (see (1)).
On the other hand, if $\pi'$ is an involution with support in $A$, then
$h\in V_{\pi'}$ and $h^2\in S_{\pi'}$ (see (2-i) of Lemma
\ref{Fremlin_Lemma1} and the definition of $V_{\pi'}$). Thus,
$\alpha(h^2)\in S_{\alpha(\pi')}$ and $supp(\alpha(h^2))\cap
supp(\alpha(\pi'))=\emptyset$ (see (4-i) of Lemma \ref{Fremlin_Lemma1}).
As $\pi'$ is arbitrary, we have from Corollary \ref{TopologyGeneration}
that $supp(\alpha(h^2))\cap \Lambda(A)=\emptyset$. Hence,
$supp(\alpha(\pi))\setminus\Lambda(A)\supset supp(\alpha(h^2))\neq
\emptyset$. $\vartriangleright$

\medskip(4) Define the map $\Lambda^{*} : RO(X_2)\rightarrow RO(X_1)$ as follows:
$$\Lambda^{*}(B)=\bigvee \{supp(\alpha^{-1}(\pi)) : \pi\in
\Gamma_2,\; \pi^2=1\mbox{ and }supp(\pi)\subseteq B\}.$$

\medskip (5) $\Lambda^{*} \Lambda (A)=A$ for every $A\in RO(X_1)$ and
$\Lambda\Lambda^{*}(B)=B$ for every $B\in RO(X_2)$.
 $\vartriangleleft$ By Lemma \ref{TopologyGeneration}, the set $A$
can be covered with clopen sets $\{C_n\}$ which are supports of
involutions from $\Gamma_1$. Take an involution $\pi$ whose support is
$C_n$. As $C_n\subset A$, $\alpha(\pi)$ is an involution with support in
$\Lambda(A)$. Hence $C_n=supp(\alpha^{-1}\alpha(\pi))$ is contained in
$\Lambda^{*}\Lambda(A)$. Since $C_n$ is arbitrary, we get that $A\subseteq
\Lambda^{*}\Lambda(A)$.

If $\pi\in \Gamma_2$ is an involution with support in $\Lambda(A)$, then
$\alpha^{-1}(\pi)$ is an involution whose support is in $A$ (see (3)).
Since $\pi$ is arbitrary, $\Lambda^{*}\Lambda(A)\subseteq A$, which shows
that $\Lambda^{*}\Lambda(A)\subseteq A$. Analogously, one can show that
$\Lambda\Lambda^{*}(B)=B$  for all $B\in RO(X_2)$. $\vartriangleright$

\medskip (6) Since $\Lambda\Lambda^{*}(A)=A$, it follows that for
any $A,B\in RO(X_1)$, $\Lambda(A)\subseteq \Lambda(B) $ iff $A\subseteq
B$. Moreover, $\Lambda$ is a bijection of $RO(X_1)$ and $RO(X_2)$.

\medskip (7) $\Lambda: RO(X_1)\rightarrow RO(X_2)$ is a Boolean algebra
isomorphism and $\Lambda^{-1}$ is its inverse. To see this, we refer the
reader to  Theorem 312L of \cite{fremlin}, which asserts that property (6)
implies that $\Lambda$ is a Boolean algebra isomorphism.

\medskip It remains only to prove that $\Lambda$ sends clopen sets onto
clopen sets.

\medskip (8) If $\pi\in \Gamma_1$ is an involution, then
$\Lambda(supp(\pi))=supp(\alpha(\pi))$. $\vartriangleleft$ Indeed, by
definition of $\Lambda$, $supp(\alpha(\pi))\subset \Lambda(supp(\pi))$. On
the other hand,
$$supp(\alpha(\pi))=\Lambda^{-1}\Lambda(supp(\alpha(\pi)))\supseteq
\Lambda(supp(\alpha^{-1}\alpha(\pi)))=\Lambda(supp(\pi)).
\vartriangleright$$

\medskip (9) $\Lambda(CO(X_1))=CO(X_2)$. $\vartriangleleft$ We note that
if $A$ is the support of an involution from $\Gamma_{1}$, then  we see
from (8) that $\Lambda(A)\in CO(X_2)$. Now if $A$ is an arbitrary clopen
set, then, by Corollary \ref{TopologyGeneration}, $A$  is a finite union
of supports of involutions. Since $\Lambda$ is a Boolean algebra
isomorphism and finite unions of clopen sets in $RO(X_2)$ coincide with
the theoretical-set unions, we obtain that $\Lambda(A)\in CO(X_2)$.
$\vartriangleright$

\medskip (10) For any $B\in RO(X_2)$ and $g\in \Gamma_1$, we have that
$\alpha(g)(B)=\Lambda g\Lambda^{-1}(B)$. $\vartriangleleft$  Assume the
converse, i.e $h=\alpha(g)^{-1}\Lambda g \Lambda^{*}$ is not the identity
automorphism of $RO(X_2)$. Notice that $h$ also preserves $CO(X_2)$. Then
there exists a clopen set $V$, the support of an involution $h\in
\Gamma_2$, such that $h(V)\cap V=\emptyset$. Let $\pi$ be an involution
from $\Gamma_2$ with support in $V$. Then $\alpha^{-1}(\pi)$ is supported
by $\Lambda^{-1}(V)$. Hence, $g\alpha(\pi)^{-1}g^{-1}$ is supported by
$g(\Lambda^{-1}(V))$ and $\alpha(g \alpha^{-1}(\pi) g^{-1})=
\alpha(g)\pi\alpha^{-1}(g)$ is supported by $\Lambda g\Lambda^{-1}(V)$. On
the other hand, $\alpha(g)\pi\alpha^{-1}(g)$ is supported by
$\alpha(g)(V)$.

Furthermore, we have that $\alpha(g)(V)\cap \Lambda
g\Lambda^{-1}(V)\neq\emptyset$. It follows that $V\cap
\alpha^{-1}(g)\Lambda g\Lambda^{-1}(V)\neq\emptyset$, which is a
contradiction. $\vartriangleright$

\medskip
Let us continue the proof of Theorem \ref{Theorem_Spatial_Realization}.
Since $\Lambda$ is an isomorphism of $CO(X_1)$ and $CO(X_2)$, it can be
extended to a homeomorphism $\widehat \alpha : X_1\rightarrow X_2$.
Moreover, it follows from (10) that $\alpha(g)(x) =\widehat\alpha
g\widehat\alpha^{-1}(x)$ for any $g\in \Gamma_1$ and $x\in X_2$. This
completes the proof of the theorem. \hfill$\square$
\medskip

To establish the fact that the group $\Gamma$ is a complete invariant of
flip conjugacy, we need to prove the following lemma. The proof we present
here is analogous to Proposition 5.8 from \cite{gps:1999}.

\begin{lemma}\label{Lemma_FullGroup} For a Cantor minimal system $(X,\f)$,
let $\Gamma$ denote any of the groups from (\ref{threeSubgroups}). Then
the topological full group of $\Gamma$ is equal to $[[\f]]$.
\end{lemma}
{\it Proof.} It is sufficient to show that $\f$ belongs to the topological
full group of $D([[\f]])$. For every $x\in X$ find a clopen neighborhood
$V_x$ such that $\f^j( V_x)\cap V_x=\emptyset$ for $j=1,2$. Define a
homeomorphism $h_x$ as follows:
$$h_x(y)=\left\{\begin{array}{ll}\f(y), & x\in V_x\cup \f(V_x)\\
 \f^{-2}(y), & x\in \f^{2}(V_x)\end{array}\right.$$
 It is not hard to see that $h_x\in D([[\f]])$ (it follows from the proof of Lemma
 \ref{ManyInvolutions}).

By  compactness of $X$, there exist $x_1,\ldots,x_n\in X$ such that
$X=\bigcup_{j=1}^n V_{x_j}$. Set $$U_1=V_{x_1},\; U_2=V_{x_2}\setminus
U_1,\ldots,\ U_n= V_{x_n}\setminus (U_1\cup\ldots\cup U_{n-1}).$$ Then
$\{U_1,\ldots, U_n\}$ forms a clopen partition of $X$. Then
$\f(x)=h_{x_i}(x)$ whenever $x\in U_i$. This completes the proof.
\hfill$\square$

\begin{theorem}\label{Theorem_FlipConjugacy}
Let $(X_1,\f_1)$ and $(X_2,\f_2)$ be Cantor minimal systems. Then $\f_1$
and $\f_2$ are flip conjugated if and only if one of the following
statements holds:

(1) $D([[\f_1]])\cong D([[\f_2]])$;

(2) $[[\f_1]]_0\cong [[\f_2]]_0$;

(3) $[[\f_1]]\cong [[\f_2]]$.
\end{theorem}

{\it Proof.} Notice that the cases (2) and (3) are already proved in
\cite{gps:1999}. Here we present a unified proof for all three cases. Let
$\Gamma_i$ denote one of the following groups: $[[\f_i]]$,  $[[\f_i]]_0$,
$D([[\f_i]])$, $i=1,2$.

It is clear that the map implementing the flip conjugacy of $\f_1$ and
$\f_2$ can be lifted up to an isomorphism between $\Gamma_1$ and
$\Gamma_2$.

Conversely, by Theorem \ref{Theorem_Spatial_Realization} every isomorphism
between $\Gamma_1$ and $\Gamma_2$ is spatial. By Lemma
\ref{Lemma_FullGroup}, this spatial isomorphism can be extended to the
 isomorphism of $[[\f_1]]$ and $[[\f_2]]$. Then it follows from  Corollary 2.7
 of  \cite{boyle_tomiyama}  that $\f_1$ and $\f_2$ are flip conjugated. \hfill$\square$

\bigskip {\it Acknowledgement}. The work was done while the authors were
visiting  Nicolas Copernicus University in Torun. The visit of
K.M. was supported by the INTAS  YSF-05-109-5315. The authors
appreciate the warm hospitality of the university and especially
thank Jan Kwiatkowski for useful discussions of the studied
problems.


\end{document}